\input{style/arxiv-ba.cfg}
\documentclass[ba,linksfromyear,preprint]{imsart}
\makeatletter
   \@ifpackageloaded{natbib}{}{\usepackage{natbib}}
\makeatother
\RequirePackage{amsthm,amsmath}
\usepackage{amssymb,algorithm,algorithmicx,algpseudocode}

\pubyear{2015}
\volume{10}
\issue{2}
\firstpage{411}
\lastpage{439}
\doi{10.1214/14-BA921}

\startlocaldefs
\setattribute{email}{text}{}

\newcommand{\logit}{ \mbox{logit}}

\newcommand{\y}{\mathbf{y}}
\newcommand{\ve}{\mathbf{e}}
\newcommand{\Sr}{\Sigma_R}
\newcommand{\hSr}{\widehat{\Sigma}_R}
\newcommand{\vS}{\mathbf{S}}
\newcommand{\vtheta}{\boldsymbol{\theta}}
\newcommand{\1}{\mathbb{I}}
\newcommand{\hf}{\widehat{f}}
\newcommand{\abcql}{ABC$_{ql}$}
\newcommand{\Ss}{\mathcal{S}}
\newcommand{\nmodel}{\pi(y\mid\vtheta)}
\newcommand{\snmodel}{\pi(y\mid\theta)}
\newcommand{\nlik}{L_{N}(\vtheta)}
\newcommand{\snlik}{L_{N}(\theta)}
\newcommand{\npost}{\pi_{N}(\vtheta\mid\y)}
\newcommand{\snpost}{\pi_{N}(\theta\mid\y)}
\newcommand{\post}{\pi^\epsilon(\vtheta\mid\mathbf{s}_{obs})}
\newcommand{\spost}{\pi^\epsilon(\theta\mid s_{obs})}
\newcommand{\prior}{\pi(\vtheta)}
\newcommand{\sqlik}{L_{Q}(\theta)}
\newcommand{\qlik}{L_{Q}(\vtheta)}
\def\Real{{\rm I\!R}}
\def\T{{{^{_{\sf T}}}}}
\endlocaldefs

\begin{document}

\begin{frontmatter}
\title{Approximate Bayesian Computation by Modelling Summary
Statistics in a~Quasi-likelihood Framework}
\runtitle{Quasi-ABC}

\begin{aug}
\author[a]{\fnms{Stefano} \snm{Cabras}\ead[label=e1]{s.cabras@unica.it}\ead[label=e4]{stefano.cabras@uc3m.es}},
\author[b]{\fnms{Maria Eugenia} \snm{Castellanos Nueda}\ead[label=e2]{maria.castellanos@urjc.es}},
\and
\author[c]{\fnms{Erlis} \snm{Ruli}\ead[label=e3]{ruli@stat.unipd.it}}

\runauthor{Stefano Cabras, Mar\'{\i}a Eugenia Castellanos, and Erlis Ruli}


\address[a]{Department of Mathematics, Universit\`{a} di Cagliari (Cagliari,
Italy), \printead{e1} and Department of Statistics, Universidad Carlos III de Madrid (Madrid,
Spain), \printead{e4}}
\address[b]{Department of Statistics and O.R., Universidad Rey
Juan Carlos (Madrid, Spain),\\ \printead{e2}}
\address[c]{Department of Statistics, Universit\`{a} di Padova
(Padova, Italy), \printead{e3}}
\end{aug}

%
\begin{abstract}
Approximate Bayesian Computation (ABC) is a useful class
of methods for Bayesian inference when the likelihood function is
computationally intractable. In practice, the basic ABC algorithm
may be inefficient in the presence of discrepancy between prior and
posterior. Therefore, more elaborate methods, such as ABC with the
Markov chain Monte Carlo algorithm (ABC-MCMC), should be used. However,
the elaboration of a proposal density for MCMC is a sensitive issue and
very difficult in the ABC setting, where the likelihood is intractable.
We discuss an automatic proposal distribution useful for ABC-MCMC algorithms.
This proposal is inspired by the theory of quasi-likelihood (QL) functions
and is obtained by modelling the distribution of the summary statistics as
a function of the parameters. Essentially, given a real-valued vector of
summary statistics, we reparametrize the model by means of a regression
function of the statistics on parameters, obtained by sampling from the
original model in a pilot-run simulation study. The QL theory is well
established for a scalar parameter, and it is shown that when the conditional
variance of the summary statistic is assumed constant, the QL has a closed-form
normal density. This idea of constructing proposal distributions is extended
to non constant variance and to real-valued parameter vectors. The
method is
illustrated by several examples and by an application to a real problem in
population genetics.
\end{abstract}

%


\end{frontmatter}


\section{Introduction}
Many statistical applications in diverse fields such as biology,
genetics and finance often involve stochastic models with analytically
or computationally intractable likeli-\break hood functions. The rapidly
growing literature on Approximate Bayesian Computation (ABC) has led to
a set of methods which do not involve direct calculation of the
likelihood, leading to Bayesian inference that is approximate in a
sense that will be specified later.

ABC methods are becoming popular in genetics \citep
{siegmund08,foll2008}, epidemiology \citep{blum2010,tanaka2006} and in
population biology \citep{ratmann2007,hamilton2005,cornuet2008} among
other areas.

Formally, let $\y=(y_1,\ldots,y_n)$ be a random sample of size $n$
drawn from a statistical model $\nmodel$ indexed by the parameter
$\vtheta\in\Theta\subseteq\Real^{p}$. The likelihood for $\vtheta
$, corresponding to $\nmodel$ is $\nlik$, which is not available in
closed expression. For a certain prior $\prior$, the aim is to obtain
the posterior distribution $\npost\propto\nlik\prior$, but as
$\nlik$ is inaccessible, $\npost$ cannot be approximated by directly
evaluating $\nlik$.

This difficulty may be overcome by using ABC methods. Specifically, let
$\mathbf{s}=s(y) \in\Ss\subseteq\Real^{p}$ be a vector of
observable summary statistics (e.g. mean, variance, quantiles etc.),
which may not be sufficient, let $\rho(\mathbf{s},\mathbf{s}_{obs})$
be a metric distance between $\mathbf{s}$ and its observed value
$\mathbf{s}_{obs}$ with $\epsilon>0$ the tolerance parameter. ABC
methods approximate $\npost$ by
\[
\post=\int_\Ss\pi^{\epsilon}(\vtheta, \mathbf{s}\mid\mathbf
{s}_{obs}) d\mathbf{s},
\]
where $\pi^{\epsilon}(\vtheta, \mathbf{s}\mid\mathbf{s}_{obs})
\propto\pi(\vtheta) \pi(\mathbf{s}\mid\vtheta) \1_{\rho
<\epsilon}$, and $\1_{\rho<\epsilon}$ is the indicator function for
the event $\left\{\mathbf{s}\in\Ss\mid\rho(\mathbf
{s}_{obs},\mathbf{s})<\epsilon\right\}$. They require a choice of
$\epsilon$ and for this purpose several authors \citep
{bortotetal2007,fasaietal2013,citeulike:12356884,
barnes2012considerate,aeschbacher2012novel} suggest approaches where
$\epsilon$ is estimated as part of an extended model with respect to
$\nmodel$, while a recent approach based on diagnostic tools for ABC
can be found in \cite{Prangle:2013uq}. In this work another criterion
for choosing $\epsilon$ is discussed.

The basic version of the ABC algorithm relies on simulation by the
mixture representation method consisting in generating, say, $T$ values
of $\vtheta$ from $\prior$ and using them to generate the
corresponding $T$ values of $\mathbf{s}$ from $\nmodel$ at the
simulated $\vtheta$. We accept all values of $\vtheta$ such that
$\rho(\mathbf{s},\mathbf{s}_{obs})<\epsilon$. For $\epsilon
\rightarrow0$ the ABC method has been proven to return a consistent
estimator of the posterior $\pi(\vtheta\mid\mathbf{s}_{obs})$ and
under some assumptions it is also possible to provide the approximation
error as shown in \cite{biau2012new}. Moreover,
if $\mathbf{s}$ is sufficient and $\epsilon\rightarrow0$, then
$\post\rightarrow\npost$. There is a certain agreement in that low
dimensional, but informative summary statistics improve the accuracy of
the ABC approximation \citep{blumetal2013}.

One drawback of the basic ABC algorithm is that it can be extremely
inefficient when the discrepancy between $\prior$ and $\npost$ is
relevant. Unfortunately, as $\nlik$ is intractable, the discrepancy
between $\prior$ and $\npost$ is difficult to know a priori and not
easy to assess. To deal with this issue, more advanced Monte Carlo
methods, such as ABC-MCMC, originally developed in \cite{Marjoram2003}
and further analyzed in several papers such as \cite
{beumontetal2009,Andrieu:2009aa,lee2012choice}, or Sequential Monte
Carlo (SMC) methods \citep[see, e.g.,][]{beaumont2009adaptive,sisson07}
may be used. All these methods attempt to account for the
observed data at the proposal stage. However, to accomplish this task,
a proposal distribution or a perturbing kernel is required, which in
practice is supplied by the analyst.\looseness=1\vadjust{\eject}

Another aspect of the method, which is of major concern in the ABC
literature, is the choice of $s$, which should be informative for
$\theta$. The same concern applies here, as we require $s$ not to be
ancillary with respect to $\theta$. Many suggestions can be found in
the current literature. For instance, \cite{fp2012} propose
considering the posterior mean of $\theta$, i.e. $E_{\post}(\theta
)$. The latter is estimated by means of a pilot-run simulation that
depends on the specific observed sample. Moreover, \cite{rulietal2013}
suggest choosing $s$ as the score of the composite likelihood function
obtained from $\nmodel$.

The present work focuses on the study of a class of proposal
distributions for ABC-MCMC, and a method for building proposal
densities, which target the posterior $\post$, is illustrated. Such
proposal distributions depend on the model at hand and account for the
observed data.
These distributions for $\theta$ are constructed in a way that leads
to adopting a normal kernel on the space of $s$, and then consider a
reparametrization from $s$ to $\theta$ by a suitable regression
function $f(\theta)=E_{\nmodel}(s\mid\theta)$. A recent approach,
strongly connected with the use of such a regression function, can be
found in \cite{citeulike:12356884}, where the $f(\theta)$ is the
binding function in indirect inference \citep{Gourierouxetal1993}.

We show that for scalar parameter problems and $f(\theta)$ such
proposal distributions arise from the class of quasi-likelihood
functions (QL) of $\theta$ \citep{mccullagh} denoted by $\sqlik$.
For multidimensional parameter problems, the QL is not tractable, but
the idea can still be generalized to these contexts using asymptotic
arguments. Indeed, for the vector of parameters $\vtheta$, we consider
a multivariate normal kernel and a multivariate transformation from
$\mathbf{s}$ to $\vtheta$. For both the scalar and the multi
parameter cases, these transformations are typically not available
analytically and we estimate them in a pilot-run simulation. This
pilot-run simulation is performed regardless of the specific observed
sample and thus it can serve for routine analysis. This is an appealing
feature of the proposed method as will be shown later by an application
to a Genome Wide Association Study (GWAS).

Despite the fact that under some more elaborate requirements for the
proposal distribution, later discussed, we end up in a proposal which
is not the QL for $\vtheta$, we think that the connection of ABC and
QL is important because $\nlik$ is not available and the estimation
theory of the QLs guarantees, asymptotically, that $\qlik$ targets
$\nlik$.

The structure of the paper is as follows: Section \ref{sec:qlabcmcmc}
discusses $\qlik$, which inspires the proposal distributions
considered throughout the paper. These distributions will be embedded
in the ABC-MCMC algorithm. The proposed ABC algorithm, \abcql, for
scalar parameters is formally discussed in Section \ref
{sec:abc-mcmc-ql}, while the generalization to $p>1$ is presented in
Section \ref{sec:abc-mcmc-ql-p}. Section \ref{sec:examples}
illustrates the proposed method with some examples from the ABC
literature and an important application to GWAS for population genetic
isolates. Conclusions and further remarks are given in Section \ref
{sec:conclusions}.

\section{The two relevant tools: quasi-likelihood and the ABC-\break MCMC
algorithm}\label{sec:qlabcmcmc}

The theory and use of estimating equations and that of the related
quasi- and quasi-profile likelihood functions have received a good deal
of attention in recent years; see, among others, \cite
{liang1995,Barndorff-Nielsen1995,desmond1997,heyde1997,adimari2002,severini2002,wang2003,jorgensen2004,bellio2008}.
In addition, \cite{a-quasi-profile,lin2006,GreRacVen08} discuss the
use of QL functions in the Bayesian setting.

Let $s=s(y) \in\Real$ be a scalar summary statistic generated from
$\pi(s(y)\mid\theta)$ whose observed value is $s_{obs}$. We assume
for convenience that the summary statistic lies on the real line, which
can be easily achieved by suitable transformations (e.g. a
log-transformation of the sample variance).

Moreover, suppose $\theta$ is a scalar parameter, i.e. $p=1$ and let
$\Psi(s;\theta)$ be an unbiased estimating function of $\theta$
based on $s$, i.e. $E_{\nmodel}\{\Psi(S;\theta)\}=0$.

The QL for $\theta$ based on $\Psi(s;\theta)$ \citep{mccullagh}, is
given by
\begin{eqnarray}
L_Q (\theta) = \exp\left\{ \int_{c_0}^\theta A(t) \Psi(s;t) \, dt
\right\}
\ , \label{quasi}
\end{eqnarray}
where $A(\theta) = M(\theta)/\Omega(\theta)$, $c_0$ is an arbitrary
constant,
\[
M(\theta) = - E\left\{\frac{\partial\Psi(S;\theta)}{\partial
\theta}\mid\theta\right\},
\]
and
\[
\Omega(\theta) = E \{\Psi(S;\theta)^{2}\mid\theta\}=\text{Var}\{
\Psi(S;\theta)\mid\theta\}.
\]
When $p = 1$, a quasi likelihood for $\theta$ is usually easy to
derive, while for $p>1$ some difficulties arise. Moreover, as shown
below, for a suitable estimating function and under $\text{Var}\{\Psi
(S;\theta)\mid\theta\}$ constant, \eqref{quasi} is a normal kernel.

The ABC-MCMC algorithm, proposed in \cite{Marjoram2003}, summarized in
Algorithm \ref{alg-mh-abc}, evaluates $\nlik$ indirectly via the
indicator function $\1_{\rho<\epsilon}$, and uses the proposal
density $q(\theta^{(t)} \mid\theta^{(t-1)})$.

\begin{algorithm}
\begin{algorithmic}[1]
\State Set $\epsilon>0$, $\theta^{(0)}=\theta_{init}\in\Theta$;
\For{$t=1$ to $T$}
\State generate $\theta^*\sim q(\theta^{(t)}|\theta^{(t-1)})$;
\State generate $s \sim\pi(s(y) \mid\theta^*)$;
\State calculate $\rho=\rho(s_{obs},s)$;
\State with probability
\[
\min\left\lbrace1,\frac{\pi(\theta^*)q(\theta^{(t-1)} \mid
\theta^{*})}{\pi(\theta^{(t-1)})q(\theta^*\mid\theta^{(t-1)})}\1
_{\rho<\epsilon} \right\rbrace
\]
accept $\theta^*$ and set $\theta^{(t)}=\theta^*$, otherwise $\theta
^{(t)}=\theta^{(t-1)}$
\EndFor
\State\Return{$\theta^{(1)},\ldots,\theta^{(T)}$}
\end{algorithmic}
\caption{The ABC-MCMC algorithm} \label{alg-mh-abc}
\end{algorithm}

Depending on how the proposal is defined, with Algorithm \ref
{alg-mh-abc} we may implement the independent Metropolis Hastings (MH)
or the Random Walk (RW) MH. SMC methods may also be considered as in
\cite{Toni:2009aa}. Finally, the proposal $q(\cdot)$ can also be
viewed as an importance function for the implementation of an
Importance Sampling (IS) simulation algorithm.

The proof that $\spost$ is the stationary distribution of Algorithm
\ref{alg-mh-abc} is contained in Theorem 1 of \cite{Marjoram2003} and
the rate of convergence depends on the choice of $q(\cdot)$ and
$\epsilon$. However, as $\snlik$ is not tractable, it is not possible
to further characterize the stochastic behavior of the induced chain.
In fact, the theoretical conditions discussed in \cite
{MengersenTweedie96} and \cite{AtchadePerron2007} may be assessed if
$\snlik$ is available in a closed form expression. Notice that the
approach can also be viewed as a pseudo marginal MH as we are working
with an estimated likelihood when evaluating the indicator function $\1
_{\rho<\epsilon}$, and results for convergence of such pseudo
marginal algorithms can also be found in \cite{Andrieu:2009aa}.

\section{The \abcql\;for a scalar parameter}
\label{sec:abc-mcmc-ql}
The main objective of this paper is to construct a proposal density
centered on the bulk of the posterior distribution. In a setting where
the likelihood cannot be computed explicitly, this issue could be
cumbersome. For this purpose, we consider a QL derived from estimating
functions based on $s$.

The following Proposition 1 provides the expression of $L_{Q}(\theta)$
for a general statistic $s$ and assuming $\text{Var}_{\snmodel}(S
\mid\theta)$ is constant.

{\bf Proposition 1}. Suppose $p=1$ and let $f(\theta)=E_{\nmodel}(S
\mid\theta)$ be a bounded regression function under the sampling
model $\nmodel$ for which the Jacobian $\mid f^\prime(\theta) \mid<
\infty$ and that the conditional variance $\text{Var}_{\nmodel}(S
\mid\theta)=\sigma^2_R$ is constant with respect to $\theta$.

Consider the following estimating function $\Psi(s_{obs};\theta
)=s_{obs}-f(\theta)$. In this case we have
\begin{equation}
L_{Q}(\theta)=\phi\left( \frac{f(\theta)-s_{obs}} {\sigma
_{R}}\right), \label{eqlq}
\end{equation}
where $\phi(\cdot)$ is the density of the standard normal distribution.

{\bf Proof}. Note that the estimating function is unbiased because $E\{
\Psi(s_{obs};\theta)\mid\theta\}= E(S \mid\theta)-f(\theta)=0$.
From the definition of $L_{Q}(\theta)$ the following quantities are needed:
\[
M(\theta)=- E\left(\frac{\partial\Psi}{\partial\theta}\mid
\theta\right)=f^\prime(\theta),
\]
\begin{eqnarray*}
\Omega(\theta) &=& E\{\Psi(s_{obs};\theta)^{2}\mid\theta\}\\
&=&\text{Var}(S-f(\theta)\mid\theta)\\
&=&\text{Var}(S \mid\theta)\\
&=&\sigma^2_R.
\end{eqnarray*}

Then $A(\theta)=f^\prime(\theta)/\sigma^2_R$ and by \eqref{quasi}
we have
\begin{eqnarray*}
L_{Q}(\theta)&=&\exp\left\{\int_{c_{0}}^{\theta} \frac{f^\prime
(t)}{\sigma^2_R}(s_{obs}-f(t))dt\right\}\\
&\propto&\frac{1}{\sigma_R}\exp\left(-\frac{(f(\theta
)-s_{obs})^{2}}{2\sigma^{2}_{R}} \right),
\end{eqnarray*}
which is the kernel of the normal distribution centered at $s_{obs}$
with variance $\sigma^2_R$. $\Box$

Expression \eqref{eqlq} suggests that $\sqlik$ is a normal density
when working in the $f(\theta)$ parametrization, and thus if one is
able to make a change-of-variable from $f(\theta)$ to $\theta$, it
could be employed as a proposal for MCMC-ABC algorithms, leading to a
broad class of ABC methods denoted \abcql\, algorithms.

Note that the constant variance assumption, $\text{Var}_{\snmodel}(S
\mid\theta)=\sigma^2_R$, leads to a closed-form proposal
distribution, a normal density. {However, } as this assumption may be
restrictive, we extend the idea of constructing proposal distributions
based on $\sqlik$, but assuming a non constant variance, $\sigma
^2_R(\theta)$, which can be estimated as well as $f(\theta)$ (see
Subsection \ref{estim.f}). The theory of QL assures that also for non
constant $\sigma^2_R(\theta)$ there exists a corresponding $\sqlik$
whose form is intractable and for this reason it cannot be used
directly as a proposal density. Instead, our proposal distribution for
a Random Walk Metropolis Hastings (RWMH) is based on a distribution of
the form of $\sqlik$ in \eqref{eqlq} where $\sigma^2_R(\theta)$ is
non constant:
\begin{equation}
q^{Q}(\theta\mid\theta^{(t-1)})= \phi\left(\frac{f(\theta
)-f(\theta^{(t-1)})} {\sigma_{R}(\theta^{(t-1)})}\right)\mid
f^\prime(\theta) \mid.\label{eq.proposal}
\end{equation}

\begin{algorithm}
\begin{algorithmic}[1]
\Require$f$, $f^\prime(\theta)$, $\sigma^2_R(\theta)$, or their
estimates ($\hf$, $\hf^\prime(\theta)$, $\widehat{\sigma
}^2_R(\theta)$).\\
Set $\epsilon>0$ and $\theta^{(0)}=f^{-1}(s_{obs})$;
\For{$t=1$ to $T$}
\State generate
\[
f^* \sim N(f(\theta^{(t-1)}),\sigma^2_R(\theta^{(t-1)}));
\]
\State set $\theta^*=\left\{\theta:f^{-1}(f^*)=\theta\right\};$
\State generate $s \sim\pi(s(y)\mid\theta^*)$;
\State calculate $\rho=\rho(s_{obs},s)$;
\State calculate the derivative, $f^\prime(\theta)$, of $f(\theta)$,
at $\theta^{(t-1)}$ and $\theta^*$;
\State with probability
\[
\min\left\lbrace1,\frac{\pi(\theta^*) q^{Q}(\theta^{(t-1)}\mid
\theta^{*}) }{\pi(\theta^{(t-1)}) q^{Q}(\theta^* \mid\theta
^{(t-1)}) }\1_{\rho<\epsilon} \right\rbrace
\]
accept $\theta^*$ and set $\theta^{(t)}=\theta^*$, otherwise $\theta
^{(t)}=\theta^{(t-1)}$
\EndFor
\State\Return{$\theta^{(1)},\ldots,\theta^{(T)}$}
\end{algorithmic}
\caption{The \abcql\, for $p=1$} \label{alg-q-mh-abc}
\end{algorithm}

Finally, this way of constructing proposals could also be extended to
other types of distributions with heavy tails, such as the $t$-Student
distribution.

On the other hand, ABC Importance Sampling (ABC-IS) can be implemented
using $q(\theta)=\sqlik|f^\prime(\theta)|$ as the importance
function from which it is possible to simulate. Assuming $\sigma
_R(\theta)=\sigma_R$, it would be enough to simulate a sample $z$
from a standard normal and then calculate $f^{-1}(z\sigma_R+s_{obs})$
to have a draw of $\theta$. The ABC-IS is completed by the evaluation
of the importance weights by means of $q(\theta)$. Also ABC-SMC can be
used starting with a sample from $q(\theta)$ and using steps S1-S3
from \cite{Toni:2009aa}. The computational requirements are almost the
same as for the ABC-IS. In fact, for ABC-SMC it is necessary to
calculate importance weights to be updated in the MC sequence.

\subsection{Estimation of $f(\theta)$, $f^\prime(\theta)$ and
$\sigma^2_R(\theta)$} \label{estim.f}
The function $f(\theta)$ can be elicited, suggested by the model in
\citet{mengersen2012approximate,ratmann2007} or by theoretical
arguments as in \cite{Heggland:2004fk}. For instance, in the genetic
model analyzed in \cite{mengersen2012approximate}, where the
constraint of the empirical likelihood plays the same role as $\Psi
(s_{obs}; \theta)$, $f(\theta)$ is built upon the score function of
the pairwise likelihood corresponding to the model. However, except in
few specific situations, $f(\theta)$, $f^\prime(\theta)$ and $\sigma
^2_R(\theta)$ are generally unknown, and we replace them in Algorithm
\ref{alg-q-mh-abc} by estimates that can be obtained in a pilot-run
simulation as stated in Algorithm \ref{alg-est-f}. In the sequel, when
referring to the \abcql\, algorithms, our intention is always that
$f(\theta)$, $f^\prime(\theta)$ and $\sigma^2_R(\theta)$ are
unknown and replaced by an estimator with the sole purpose of providing
the input as a proposal density for Algorithm \ref{alg-q-mh-abc} (or
Algorithm \ref{alg-q-mh-abc-multi}). If $f(\theta)$, $f^\prime
(\theta)$ and $\sigma^2_R(\theta)$ were known, the computational
effort for the \abcql\, algorithm would be reduced. In any case, the
proof of the convergence of Algorithm \ref{alg-q-mh-abc} (or Algorithm
\ref{alg-q-mh-abc-multi}) is discussed in the previous Section \ref
{sec:qlabcmcmc}.

\begin{algorithm}
\begin{algorithmic}[1]
\Require$M$, $\tilde\Theta$
\State consider $M$ values $\tilde\theta=(\tilde\theta_1,\ldots
,\tilde\theta_M)$ taken in a regular spaced grid of a suitable large
subset $\tilde\Theta\subseteq\Theta$;
\State generate $\tilde s=(\tilde s_1,\ldots,\tilde s_M)$ where
$\tilde s_m\sim\pi(s(y)\mid\tilde\theta_{m})$;
\State regress $\tilde s$ on $\tilde\theta$ obtaining $\hf(\theta)$
and $\hf^{\prime}(\theta)$;
\State regress $\left\{\log\left( \hf(\tilde\theta_m)-\tilde s_m
\right)^2\right\}_{m=1,\ldots,M}$ on $\tilde\theta$ obtaining
$\widehat{\sigma}^2_R(\theta)$.
\State\Return$\hf(\theta)$, $\hf^{\prime}(\theta)$ and $\widehat
{\sigma}^2_R(\theta)$.
\end{algorithmic}
\caption{Estimation of $f(\theta)$, $f^\prime(\theta)$ and $\sigma
^2_R(\theta)$ for $p=1$} \label{alg-est-f}
\end{algorithm}

The functions $\hf(\theta)$ or $\widehat{\sigma}^2_R(\theta)$ can
be any estimator which provides smoothing regression functions, and
$\hf(\theta)$ is at least once differentiable. This implies that the
main assumption for $f(\theta)$ is to be monotone and once
differentiable. We find it useful to use smoothing splines, for which
the derivative, $\hf^\prime(\theta)$, can be obtained analytically
from splines coefficients. Other choices are possible and left to the
convenience of the analyst that inspects the scatter diagram of points
$\left\{\tilde s_m,\tilde\theta_m\right\}_1^M$ and provides a
goodness-of-fit argument that justifies the choice. The inverse $\hf
^{-1}(f^*)$, at some point $f^*$, can be obtained either analytically,
with the bisection method on $\hf(\theta)=f^*$ or by numerical
minimization of $(\hf(\theta)-f^*)^2$, e.g., by a Newton-Raphson algorithm.

Some observations are appropriate.
\begin{itemize}
\item[$i$)] Since we are able to simulate from $\pi(y\mid\theta)$,
then $\hf(\theta)$ and $\widehat{\sigma}^2_R(\theta)$ can be
practically estimated with a precision that depends on the available
computational resources. Also more precision can be achieved by making
the regular grid $\tilde\Theta$ wider or by increasing the number of
simulations, $M$. Values of $M$ ranging from $100$ to $1000$ are enough
in the examples discussed later for $p=1$, while larger values are
needed for $p>1$, due to the curse of dimensionality as the
computational effort increases exponentially with $p$.
\item[$ii$)] The range of $\tilde\Theta$ should always be large
enough to include the observed $s_{obs}$ in order to gain precision
around the bulk of the target posterior $\spost$.
\item[$iii$)] The monotonicity assumption is a necessary condition for
ABC as it states that there exists a relation between $s$ and $\theta$
through $f$, see e.g., \cite{ratmann2007}. The lack of monotonicity is
not a fault of the proposed method, but instead it is an indication of
the fact that $s$ is not informative for $\theta$ in some subset of
the parameter space. This would be automatically recognized by the
proposal as it would be essentially flat in such a region due to a
small Jacobian.
\item[$iv$)] The function $\sqlik$ can also be useful to fix
$\epsilon$, because under the assumption that $\snpost$ is $\snlik$
dominated and that $\sqlik$ was its approximation, then it would be
enough to simulate $\theta$ from $\sqlik$ and $s$ from $S(y)|\theta$
obtaining thus the distribution for $\rho(s_{obs},s)$ and fixing
$\epsilon$ as its suitable quantile.
\item[$v$)] The computational cost for approximating $\npost$ with
the proposed approach should include that for estimating $\hf(\theta
)$, approximating the inverse $\hf^{-1}(f^*)$ and its Jacobian $\hf
^{\prime}$, where the latter is mainly important for $p>1$, as for
$p=1$ the derivative is obtained analytically. Such costs in the
implemented examples are actually reduced by another interpolation of
the inverse and Jacobian with splines, or its corresponding equivalent
Generalized Additive Model (GAM) for $p>1$ (see the next section). This
interpolation speeds up the MCMC because at each step we do not need to
calculate its inverse and Jacobian, but just its interpolation.
Finally, because the most important cost is model simulations rather
than regression estimation, we note that, in the analyzed examples, the
number of simulated statistics $M$ is no larger than 10\% of the number
of MCMC steps $T$. Such computational effort is enough to achieve the
desired precision in formulating a proposal distribution.
\end{itemize}

\section{The \abcql\;for $p>1$}\label{sec:abc-mcmc-ql-p}
Suppose $p>1$ and let $\mathbf{s}_{obs}$ be the vector of the $p$
observed statistics. In this case the theory for $\qlik$ is not well
developed for finite samples; however, when $p > 1$, $\qlik$ exists if
and only if the matrix $M(\vtheta)$ is symmetric. In the multi
parameter case, we use the regular asymptotic argument for the
likelihood \citep[see e.g.,][Ch.~4]{PaceSalv97:PrincStatiInfer}. That
is, for $n \rightarrow\infty$, the Taylor expansion of the QL around
its mode leads to the following multivariate normal QL:
\[
L_Q(\vtheta) = \Sr^{-\frac{1}{2}}\exp\left(-\frac{1}{2}
(f(\vtheta)-\mathbf{s}_{obs})\T\Sr^{-1}(f(\vtheta)-\mathbf
{s}_{obs}) \right),
\]
where $f(\vtheta)=E(\vS\mid\vtheta)$ is a bounded monotone, and
possibly non-linear regression function and $\Sr$ is the conditional
covariance matrix of $\vS\mid\vtheta$. Following the same approach
as for $p=1$, we use $L_Q(\vtheta)$ as a proposal distribution to be
used in a MH scheme (see Algorithm \ref{alg-q-mh-abc-multi}) also
considering a non constant covariance matrix $\Sr(\vtheta)$, that is
\begin{equation}
q^{Q}(\vtheta\mid\vtheta^{(t-1)})= N_p\left(f(\vtheta^{(t-1)}),
\Sr(\vtheta^{(t-1)})\right)\mid J(\vtheta) \mid, \label{eq.proposalmulti}
\end{equation}
where $N_p(\cdot,\cdot)$ denotes the $p$-variate normal distribution
with mean $f(\vtheta^{(t-1)})$ and vari\-ance-covariance matrix $\Sr
(\vtheta^{(t-1)})$, and $J(\vtheta)$ is the Jacobian of the
transformation $f(\vtheta)=E(\vS\mid\vtheta)$.

In the multi-parameter case we have a system of non-linear equations
$f(\vtheta)=\mathbf{s}$ whose solution is $\vtheta=f^{-1}(\mathbf
{s})$. Such a solution and the calculation of the determinant of the
Jacobian, $|J(\vtheta)|$, can be obtained using numerical methods for
solving a non-linear system of equations and approximating the
derivative of $f(\vtheta)$ at $\vtheta$. In the case of non constant
covariance matrix, we use the proposal distribution with covariance
\[
\Sr(\vtheta)=\left(
\begin{array}{ccc}
\sigma^2_{R1}(\vtheta) & 0 & 0\\
0 & \ddots& \\
0 & 0 & \sigma^2_{Rp}(\vtheta) \\
\end{array}
\right) ,
\]
where $\sigma^2_{R1}(\vtheta) ,\ldots,\sigma^2_{Rp}(\vtheta)$ are
the conditional variance functions for each component of $\mathbf{s}$
with respect to all $p$ components of $\vtheta$. Note that in this
case we are forced to use a diagonal covariance matrix in order to
guarantee that $\Sr(\vtheta)$ is positive definite. The correlation
between the $p$ parameters is then accounted for in the MCMC sampling.
Algorithm \ref{alg-est-fmulti} illustrates how to obtain estimates of
$f(\vtheta)$ and $\Sr(\vtheta)$ along with the calculation of the
Jacobian corresponding to $\hf(\vtheta)$ for $p>1$.

\begin{algorithm}
\begin{algorithmic}[1]
\State Consider the set of $m=1,\ldots,M^p$ points $\tilde\vtheta
_{m}=(\tilde\theta_{1m}, \ldots, \tilde\theta_{pm})$, each of $p$
scalar coordinates over a regular lattice of $\Theta_1 \times\ldots
\times\Theta_p$ and let $\tilde\vtheta$ be the $M^p \times p$
matrix of all points;
\State Generate $\mathbf{\tilde{s}}_{m} \sim\pi(\vS(y)\mid\tilde
\vtheta_m)$ and let $\mathbf{\tilde{s}}$ be the $M^p \times p$ matrix
of all simulated statistics;
\ForAll{$j=1,\ldots,p$} Regress column $j$ of $\mathbf{\tilde{s}}$,
$\tilde{s}_j$, on $\tilde\vtheta_{m}$ obtaining $\hf_j(\vtheta)$ and
regression residuals $\ve_j$. Calculate $\hat{J}(\vtheta)$ using
Richardson's extrapolation (using R package {\tt numDeriv}). \EndFor
\State Let $\ve=(\ve_1,\ldots,\ve_p)$ be the $M^p \times p$ matrix
of regression residuals,
\If{$\Sr(\vtheta)=\Sr$ constant} \State calculate $\hSr=M^{-1}\ve
\T\ve$; \EndIf
\If{$\Sr(\vtheta)$ is non constant} \State regress $\log(\ve_j^2)$
on $\tilde\vtheta_{m}$ to have $\hat{\sigma}^2_{Rj}(\vtheta)$ for
$j=1,\ldots,p$ and obtain $\hSr(\vtheta)$. \EndIf
\State\Return$\hf(\vtheta)=(\hf_1(\vtheta_1),\ldots,\hf
_p(\vtheta_p))$, $\hat{J}(\vtheta)$ and $\hSr(\vtheta)$.
\end{algorithmic}
\caption{Estimation of $f(\vtheta)$, $J(\vtheta)$ and $\Sr(\vtheta
)$ for $p>1$} \label{alg-est-fmulti}
\end{algorithm}

Once we have estimated $f(\vtheta)$ and $\Sr$ with $\hf(\vtheta)$
and $\hSr$, respectively, we can calculate the proposal $q^Q(\cdot)$
in \eqref{eq.proposalmulti} and apply Algorithm \ref
{alg-q-mh-abc-multi}. This is just Algorithm \ref{alg-q-mh-abc} in its
multivariate version, where the distance function $\rho: \Real^p
\rightarrow\Real^+$ must consider the joint distance of all $p$
coordinates of $\mathbf{s}$ with respect to $\mathbf{s}_{obs}$. Also
here $\qlik$ can be used to fix $\epsilon$ in two ways. The first
solution, which is the one adopted in this paper, consists of
considering a common $\epsilon$ for all $p$ dimensions by
characterizing the stochastic norm of $\|\mathbf{s}-\mathbf{s}_{obs}\|
$ and its quantiles. The second solution would be to consider different
tolerance parameters, one for each of the $p$ dimensions, by deriving
the $p$ marginals from the joint proposal and then, at each iteration
$t$, updating each one of the $p$ parameters separately.

\begin{algorithm}
\begin{algorithmic}[1]
\Require$f$, $J(\theta)$ and $\Sr(\theta)$, or their estimates
($\hf$, $\hat{J}(\theta)$ and $\hSr(\theta)$).
\State Set $\epsilon>0$, $\vtheta_{0}=f^{-1}(\mathbf{s}_{obs})$;
\For{$t=1$ to $T$}
\State generate
\[
f^* \sim N_p\left(f(\vtheta^{(t-1)}),\Sr(\vtheta^{(t-1)})\right);
\]
\State set $\vtheta^*=\left\{\vtheta:f^{-1}(f^*)=\vtheta\right\};$
\State generate $\mathbf{s}\sim\pi(\mathbf{s}(y)\mid\vtheta^*)$;
\State calculate $\rho=\rho(\mathbf{s}_{obs},\mathbf{s})$;
\State calculate the determinant of the Jacobian matrices $J(\vtheta
^{(t-1)})$ and $J(\vtheta^{*})$;
\State with probability
\[
\min\left\lbrace1,\frac{\pi(\vtheta^*) q^{Q}(\vtheta^{(t-1)}
\mid\vtheta^*) }{\pi(\vtheta^{(t-1)})q^{Q}(\vtheta^* \mid\vtheta
^{(t-1)}) }\1_{\rho<\epsilon} \right\rbrace
\]
accept $\vtheta^*$ and set $\vtheta^{(t)}=\vtheta^*$, otherwise
$\vtheta^{(t)}=\vtheta^{(t-1)}$
\EndFor
\State\Return{$\vtheta^{(1)},\ldots,\vtheta^{(T)}$}
\end{algorithmic}
\caption{The \abcql\; for $p>1$} \label{alg-q-mh-abc-multi}
\end{algorithm}

Finally, in order to speed up the MCMC algorithm, especially for large
$p$, it is worth noting that once $\hf(\vtheta)$ is estimated, its
inverse and the Jacobian $\hat{J}(\vtheta)$ can be further
interpolated by means of their respective values calculated on the
points of $\vtheta$ used for the pilot-run. With such interpolation,
the inverse and Jacobians are calculated only on the points of the grid
($M^p$ in total), which is much less computationally demanding than a
calculation for all MCMC steps.

Many of the remarks outlined in $i$)-$v$) hold in the multivariate case
as well. In \mbox{order} to guarantee enough flexibility, and because we are
mainly interested in predicting $\mathbf{s}$, we consider $\hf
(\vtheta)$ and $\hat{\sigma}^2_{R1}(\vtheta),\ldots,\hat{\sigma
}^2_{Rp}(\vtheta)$, to belong to the class of generalized additive
regression models \citep{stone1985} in which each component of
$\vtheta$ enters into the linear predictor by means of a smoothing
spline as discussed, for instance, in Section 12.2 of \cite
{faraway2006}. The Jacobian of the non linear system, $\hat{J}(\vtheta
)$ which relates $p$ summary statistics to the $p$ parameters, is
calculated using Richardson's extrapolation (implemented in the R
package {\tt numDeriv}). Finally, the inverse of the non linear system
of equations at some point $\mathbf{s}$ is obtained by Newton steps as
detailed in \cite{dennis1996} (implemented in the R package {\tt nleqslv}).

We acknowledge that the approach proposed here limits the number of
observed statistics to be equal to the number of unknown parameters.
This is in line with the general recommendation to keep the number of
statistics, i.e. the number of estimating functions, equal to the
number of parameters, as also discussed in \cite
{mengersen2012approximate} and \cite{rulietal2013}.

\section{Examples}\label{sec:examples}
In this section we illustrate the proposed approach with four examples.
The first is a coalescent model \citep{Tavare:1997fk} with a scalar
parameter of interest. The second example is a gamma model with two
unknown parameters, and the third is the $g$-and-$k$ distribution
\citep[see, e.g.,][]{mcvinish2012} with four unknown parameters. The
last example is an application to a real dataset concerning GWAS, with
three unknown parameters. In the first example we apply Algorithms \ref
{alg-q-mh-abc} and \ref{alg-est-f}, whereas for the other examples we
apply Algorithms \ref{alg-est-fmulti} and \ref{alg-q-mh-abc-multi}.

The coalescent model and four-parameter $g$-and-$k$ distribution are
considered as benchmark examples. The example of the gamma model is
useful in order to validate the procedure against a known $\nlik$.
Finally, the example of GWAS is relevant for analyzing population
genetic isolates for which the genealogy tree is known. Although these
kinds of data are rare, they are exceedingly more powerful for
detecting genes related to some phenotypes than the usual and more
costly GWAS applied to larger samples of open populations. In all the
examples we have used RW-ABC-MCMC, although in Examples 1 to 3 we have
used also the independent MH ABC algorithm, obtaining similar results.
This is because the proposal distribution, constructed under the
framework of QL functions, is located around the bulk of the posterior
distribution, $\post$. Finally, the use of proper priors, in some
examples, has the sole purpose of allowing comparison with existing
methods. All examples, except the GWAS, have been implemented assuming
that $\widehat{\sigma}^2_R(\theta)$ and $\hSr(\vtheta)$ are non
constant. In the sequel \abcql\,denotes the types of algorithms where
the proposal for the MCMC has been approximated as explained above.

\subsection{Coalescent model}
In the following we consider an example from population genetics,
namely the coalescent model analyzed in \cite{Tavare:1997fk} and \cite
{blum-francois}, among others. Given a set of $n$ DNA sequences, the
aim is to estimate the effective mutation rate, $\theta^\prime>0$,
under the infinitely-many-sites model. In this model, mutations occur
at rate $\theta^\prime$ at DNA sites that have not been hit by
mutation before. If a site is affected by a mutation, then it is said
to be segregating in the sample. In this example, the summary statistic
$s^\prime=y$ is the number of segregating sites. The generating
mechanism for $s^\prime$ is the following:
\begin{itemize}
\item[1] Generate $T_{n}$, the length of the genealogical tree of the
$n$ sequences, where $T_{n}=\sum_{j=2}^{n} j W_{j}$, where $W_{j}$ are
independent Exponential random variables with mean $2/j(j-1)$ such that
$T_n$ has mean $\mu_{T_n}=2\sum_{j=1}^{n-1}1/j$ and variance \\
$\sigma_{T_n}^2=4\sum_{j=1}^{n-1}1/j^2$;
\item[2] Generate $(S^\prime\mid\theta^\prime,T_n) \sim
Poisson(\theta^\prime T_{n}/2).$
\end{itemize}
Hence, the likelihood $L_N(\theta^\prime)$ is given by the marginal
density of $(S^\prime\mid\theta^\prime)$ with respect to $T_n$,
which has a closed form only for $n=2$ as $T_{2}\sim Exp(1/2)$. For
large $n$, we approximate the inner integral in $T_n$ by simulating
$10^5$ $W_j$ for $j=1,\ldots,n$ and then obtaining the marginal
density of $(S^\prime\mid\theta^\prime)$ by averaging over these
$10^5$ simulated values of $T_n$. This parametric approximation is
denoted by $\pi_{ap}(\theta^\prime\mid s^\prime_{obs})$ and relies
on the partial knowledge of the likelihood as the marginal density is
obtained using the Poisson likelihood.

\begin{figure}
\includegraphics{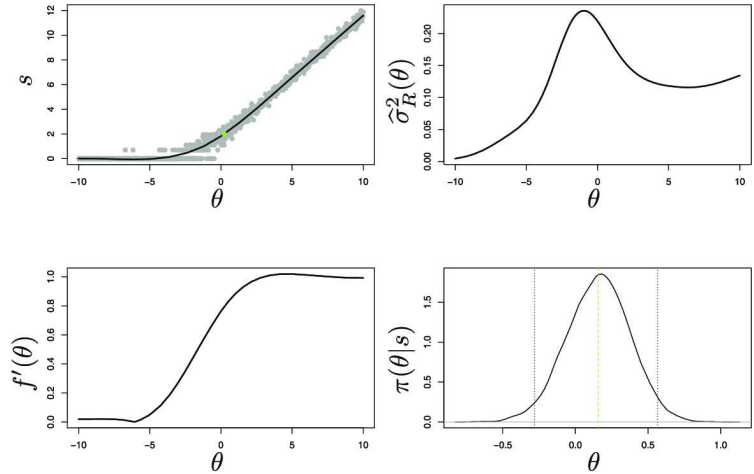}
\caption{Example: coalescent model. (Top-left) Realizations of $s$ for
the pilot-run study along with the estimated $\hf$ and $s_{obs}= 2$
(green). (Top-right) The estimated conditional variance of $s|\theta$,
$\widehat{\sigma}^2_R(\theta)$. (Bottom-left) The Jacobian of $\hf
(\theta)$ in the grid of the pilot-run. (Bottom-right) The
approximated posterior $\pi(\theta|s_{obs})$ with 95\% credible
interval and posterior mean (green).}
\label{fig-coal0}
\end{figure}

In order to employ \abcql, we consider $\theta=\log(\theta^\prime
)$ and $s=\log(s^\prime+1)$ instead of $\theta^\prime$ and
$s^\prime$, respectively. As an example, we consider $n=100$ and an
observed value of $s_{obs} = 2$, which is a value likely to be obtained
under $\theta=0$ (and hence $\theta^\prime=1$). Figure \ref
{fig-coal0} reports the result of the pilot-run study (top-left) with
$M=1000$ with $\hf$, $\widehat{\sigma}^2_R(\theta)$ (top-right),
the calculated Jacobian with the splines coefficients (bottom-left) and
the approximated posterior (bottom-right) with $\epsilon$ being the
10\% quantile of the distribution of $\rho(s,s_{obs})$ simulated from
$S | \theta$, where $\theta\sim\sqlik$.

From Figure \ref{fig-coal0} we can see that the chosen summary
statistic is not informative for small values of $\theta$ because
observing no segregating sites with $n=100$ samples may occur for
almost every mutation rate lower than $e^{-5}$. This is, of course, not
a fault of the method, but of the chosen summary statistic and it may
also occur in the standard original ABC approaches. Moreover, this is
reflected by the proposed method as the estimated Jacobian is near 0
for values of $\theta^\prime<e^{-5}$. It can also be seen that for
larger values, the approximated Jacobian is nearly constant, suggesting
that there exists a linear relation between $s$ and $\theta$.

\begin{figure}
\includegraphics{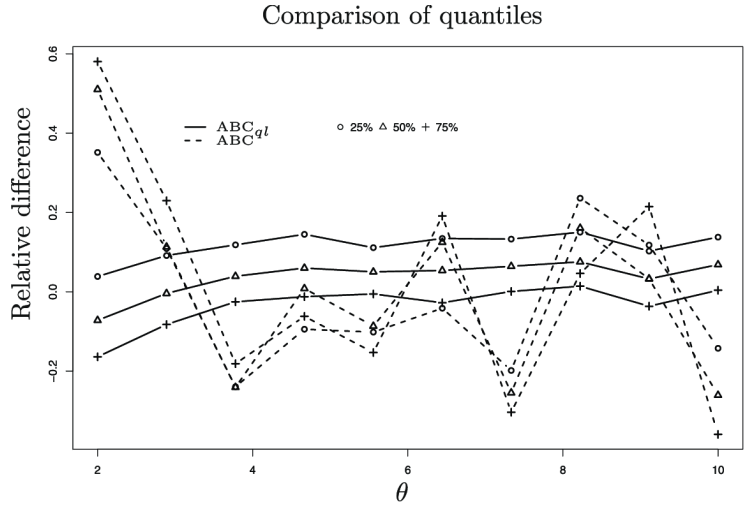}
\caption{Example: coalescent model. Comparison of ABC and \abcql\,
with the parametric approximation in terms of relative differences
between quantiles of the parametric approximation and those of ABC and
\abcql.}
\label{fig-coal2}
\end{figure}

For $\pi(\theta)=Exp(1)$ and $n=100$ we calculated, for each dataset
simulated at different values of $\theta\in(2,3,\ldots,10)$, the
relative difference of quantiles of each posterior with respect to
those obtained with the parametric approximation, as in \cite
{blum-francois}. The relative difference is defined as
$(Q_p-Q^0_p)/Q^0_p$ where $Q_p$ and $Q_p^0$ are the $p$-th quantiles of
the ABC posterior and that of the parametric approximation,
respectively. Figure \ref{fig-coal2} shows the relative differences.
We can clearly see that these differences are more robust with respect
to $\theta$ for the \abcql\,rather than for the ABC, and this is due
to the impact of the prior in the standard ABC algorithm. In fact, for
$\theta\rightarrow\infty$ the discrepancy between prior and
posterior becomes important.

\subsection{Gamma model with unknown shape and scale parameters}
Let $y \sim Gamma(\theta_1,\theta_2)$ with mean $\exp(\theta
_1-\theta_2)$ and variance $\exp(\theta_1-2\theta_2)$, where
$\theta_1$ and $1/\theta_2$ are the log of shape and scale,
respectively. We have $p=2$, $\vtheta\in\Theta= \Theta_{1} \times
\Theta_{2} \equiv\Real^{2}$ and consider the following two
statistics $\mathbf{s}=(s_1,s_{2})$, where $s_1$ and $s_{2}$ are the
logarithms of the sample mean and the standard deviation, respectively.
For $\vtheta=(0,0)$ we consider a sample of size $n=10$ with $\mathbf
{s}_{obs}=(-0.12,-0.26)$ and estimate the posterior distribution under
two independent standard normal priors, $\pi(\vtheta)=\pi(\theta
_1)\times\pi(\theta_2)=N(0,1)\times N(0,1)$. We consider the
estimation of $\hf(\vtheta)$ with $M=100$, $p=2$ over a regular
lattice of $M^2=10^4$ points in $(-2,2) \times(-2,2) \in\Theta_{1}
\times\Theta_{2}$.

\begin{figure}[t]
\includegraphics{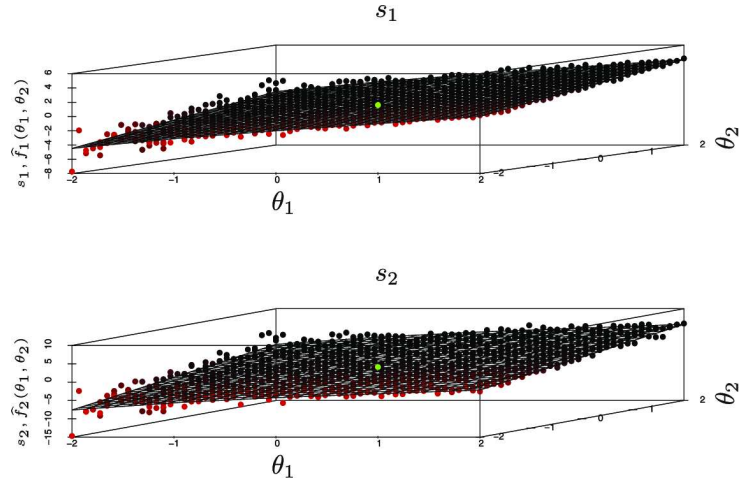}
\caption{Example: gamma model. Conditional regression functions of
each statistic against the two parameters. The green point is the
observed value.}
\label{fig-gamma1}
\end{figure}

\begin{figure}[t!]
\includegraphics{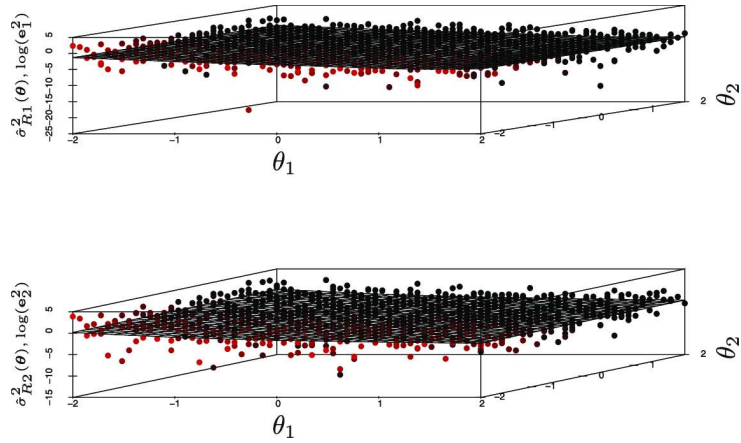}
\caption{Example: gamma model. Conditional variance functions of each
statistic against the two parameters.}
\label{fig-gamma-varest}
\end{figure}

\begin{figure}
\includegraphics{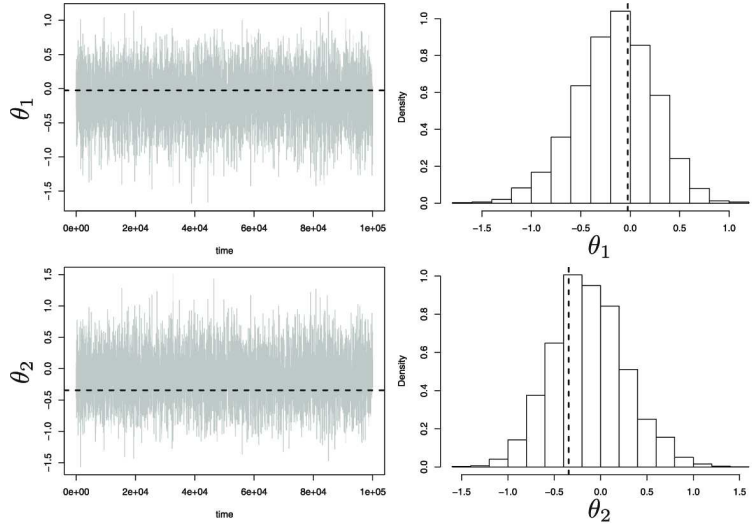}
\caption{Example: gamma model. Marginal output of Algorithm \ref
{alg-q-mh-abc-multi}, for $\theta_1$ and $\theta_2$ along with the
histogram of the marginal posteriors $\pi^\epsilon(\theta_1\mid
s_{obs})$ and $\pi^\epsilon(\theta_2\mid s_{obs})$. Vertical dotted
lines are the true values of $\theta_1$ and $\theta_2$ that are used
to generate $s_{obs}$.}
\label{fig-gamma3}
\end{figure}

\begin{figure}
\includegraphics{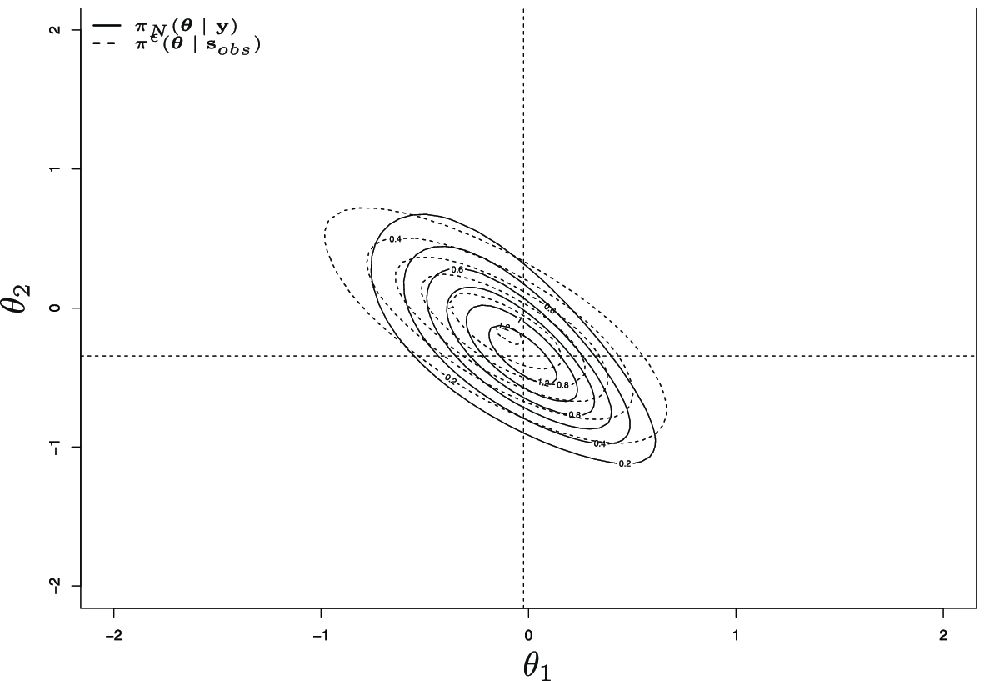}
\caption{Example: gamma model. Contours for $\post$ (dashed) along
with $\npost$ (continuous). The posterior modes of $\theta_1$ and
$\theta_2$ are represented by the cross point of the dashed lines,
where the true value is $\vtheta=(0,0)$.}
\label{fig-gamma4}
\end{figure}

Figure \ref{fig-gamma1} shows the conditional regression functions of
each statistic against the parameters, which appear quite linear.
Figure \ref{fig-gamma-varest} illustrates the $\log$ of squared
residuals and the estimation of the conditional variances $\hat{\sigma
}^2_{R1}(\vtheta)$ and $\hat{\sigma}^2_{R2}(\vtheta)$ with respect
to $\vtheta$. Figures \ref{fig-gamma3} and \ref{fig-gamma4}
illustrate the MCMC output. From the former we can deduce that the
chain mixes well, and that marginal posteriors $\pi^\epsilon(\theta
_1\mid s_{obs})$ and $\pi^\epsilon(\theta_2 \mid s_{obs})$ are
centered around the true values. From the latter, we see that the
bivariate density $\post$ obtained by the \abcql, and the true
underlying posterior $\pi_N(\vtheta\mid y)$ are similar in terms of
contour levels.

We notice that the output of the MCMC in Figure \ref{fig-gamma4} is
similar to that obtained using $\Sr$ constant (not reported here).
This is presumably because, in this example, the use of logarithms on
the scale of the summary statistics stabilizes their conditional
variance with respect to $\vtheta$. Another reason is because the
chain moves in a radius of $\mathbf{s}_{obs} < \epsilon$ where the
variance functions $\hat{\sigma}^2_{R1}(\vtheta)$, $\hat{\sigma
}^2_{R2}(\vtheta)$ are almost constant and do not differ significantly
from the estimated diagonal of $\Sr$ when it is assumed constant.

\subsection{Four-parameter $g$-and-$k$ distribution}
Distributions based on quantiles are of great interest because of their
flexibility. However, although a stochastic representation is
available, their density and hence the likelihood $L_N(\vtheta)$ are
not available in closed form and, in general, they are difficult to evaluate.

We focus on the four-parameter ($p=4$), $g$-and-$k$ distribution, which
has the following stochastic representation
\begin{eqnarray*}
z& \sim& N(0,1)\\
(y\mid z,\theta_1,\theta_2,\theta_3,\theta_4) &=& \theta_1 + \exp
{(\theta_2)} \left( 1 + 0.8 \frac{1-\exp(-\theta_3z)}{1+\exp
(-\theta_3z)}\right)\left(1+z^2 \right)^{\exp{(\theta_4)-1/2}}.
\end{eqnarray*}
The unknown parameters $\theta_1, \exp{(\theta_2)}, \theta_3$ and
$\exp{(\theta_4)}-1/2$ represent location, scale, skewness and
kurtosis, respectively \citep{haynes1997}.

Such distributions have also been used for testing several ABC
approaches as in \cite{Marjoram2003,mcvinish2012}.

We consider the following four statistics all based on empirical
quantiles $q$ of $y_q$:
\[
s=\left(y_{0.5},\log(y_{0.75}-y_{0.25}),\frac
{y_{0.75}+y_{0.25}-2y_{0.5}}{y_{0.75}-y_{0.25}}, \log\left
(y_{0.975}-y_{0.025}\right)\right)\in\Real^4,
\]
with the following meaning: $s_1$ is the median, $s_2$ is the log of
the interquartile range and $s_3$ is the skewness index described in
\cite{Bowley:1937fk}, a special case of the \cite{hinkley1975} index.
Finally $s_4$ is the log transformation of the kurtosis index described
in \cite{crowl1967}.

\begin{figure}[t]
\includegraphics{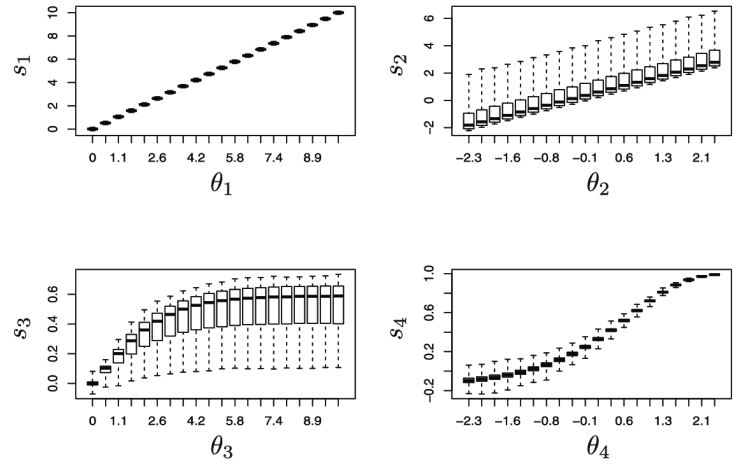}
\caption{Example: $g$-and-$k$ distribution. Conditional distribution
of $s_i|\theta_i$, for $i=1,\ldots,4$. Each boxplot represents the
marginal distribution of the statistic at the specified value of
$\theta_i$ in the horizontal axis.}
\label{fig-gk1}
\end{figure}

\begin{figure}[t!]
\includegraphics{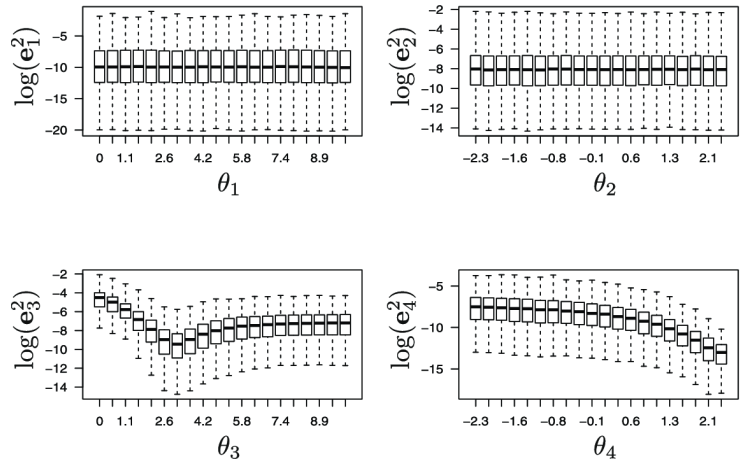}
\caption{Example: $g$-and-$k$ distribution. Logarithms of squared
residuals in the pilot-run. Each boxplot is the conditional
distribution of $\log(\ve^2_j)|\theta_j$, $j=1,2,3,4$. Such values
indicate that $\Sr(\vtheta)$ is not constant with respect to $\vtheta$.}
\label{fig-gk1-var}\vspace*{-4pt}
\end{figure}

From the pilot-run simulation shown in Figure \ref{fig-gk1}, it is
possible to notice that the relationship between statistics and
parameters is not linear. In this example, we consider a sample
simulated under the scenario discussed in \cite{fp2012} in which
$n=10^4$ observations are generated with $\theta=(3,0,2,-\log(2))$
and the uniform prior on $[0,10] \times[-\log(10),\log(10)] \times
[0,10] \times[-\log(10),\log(10)]$ is considered. Figure~\ref
{fig-gk1-var} reports the conditional distributions of the logarithm of
squared residuals, indicating that the variance may be non constant.
Therefore, for this model, we estimated with GAMs the conditional
variance functions $\hat{\sigma}^2_{Rj}(\vtheta)$, for $j=1,\ldots
,4$ as explained in Algorithm \ref{alg-est-fmulti}.

\begin{figure}[t]
\includegraphics{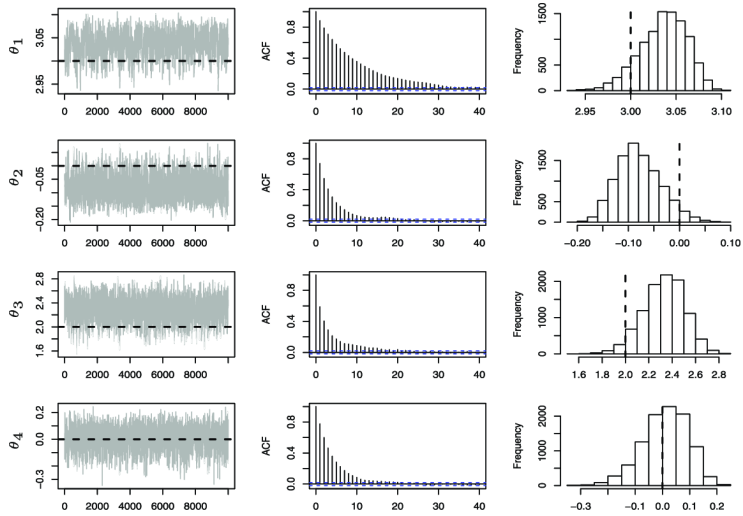}
\caption{Example: $g$-and-$k$ distribution. The approximated marginal
posterior $\pi^\epsilon(\theta_i\mid\mathbf{s}_{obs})$, $i=1,2,3,4$.}
\label{fig-gk2}
\end{figure}

Figure \ref{fig-gk2} shows the output of the \abcql\,with the
RW-ABC-MCMC algorithm for the four marginal posterior densities, for a
given sample. We can see that the marginal posterior distributions
include, in their high posterior density interval, the true value of
$\vtheta$ and the chains have good mixing.

\begin{figure}[t!]
\includegraphics{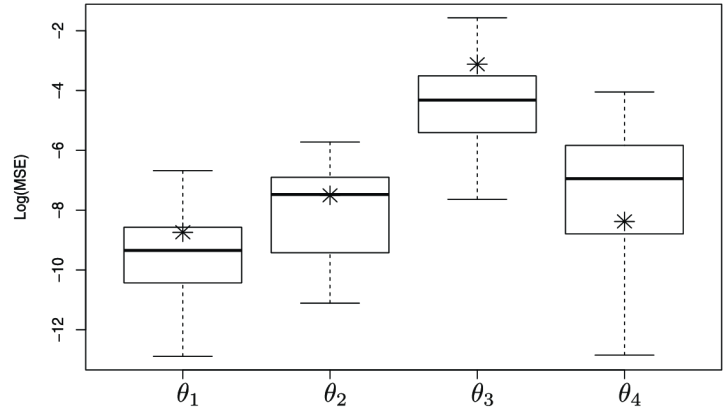}
\caption{Example: $g$-and-$k$ distribution. The marginal distributions
of the Log of MSE for each of the four parameters, in 50 replications
of the considered simulation scenario, along with the Log of MSE of the
Semi-Automatic ABC procedure reported in \cite{fp2012} (stars).}
\label{fig-gk3}
\end{figure}

Finally, the performance of the proposed method is assessed by a
simulation study~of 50 simulated datasets. As in \cite{fp2012}, the 50
datasets are simulated from 50 different parameter values sampled from
the prior. The expected quadratic errors for each marginal \abcql$\,$
posterior are reported in Figure \ref{fig-gk3}, where also the Mean
Squared Error (MSE) of the Semi-Automatic ABC, taken from \cite
{fp2012}, is shown for comparison of the order of magnitudes.

From Figure \ref{fig-gk3}, we can see that the expected quadratic
error under \abcql\; is compatible with that reported for the
Semi-Automatic ABC. Notice that our method uses a set of four
observable summary statistics that differ from those ones used in
\cite{fp2012}.

\subsection{GWAS for isolated populations with known genealogy}
In this application, we address the problem of estimating DNA markers
related to a certain phenotype such as, for instance, the presence of a
certain disease. In this problem genotype is represented by a large set
of DNA sequences known as Single-Nucleotide Polymorphisms or SNPs in
the sequel. Such SNPs are usually observed in millions per individuals
and thus fast and reliable statistical methods are needed in order to
answer the scientific question as to which SNPs are mainly related to
the disease. A~dataset for a case/control study is usually collected on
an open population where the degree of inbreeding, that is the mating
of pairs who are closely related genetically, is unknown as it is
usually negligible in open populations. Nonetheless, there exists
certain evidence in the genetic literature that very valuable
indications regarding SNPs/Disease\vadjust{\eject} relationships may come from the
study of isolate genetics, i.e. human samples for which inbreeding is
also relevant and known. Such types of collected samples are very rare
because there are very few genetic isolates in the world, and so the
statistical methods to analyze them are not very well developed. One
example of a genetic isolate for which data are also available is the
Sardinian genetic isolate of the Ogliastra region, which is situated in
the center of the island of Sardinia \citep{a-cabras-etal-bmc}.

\begin{figure}[t]
\includegraphics{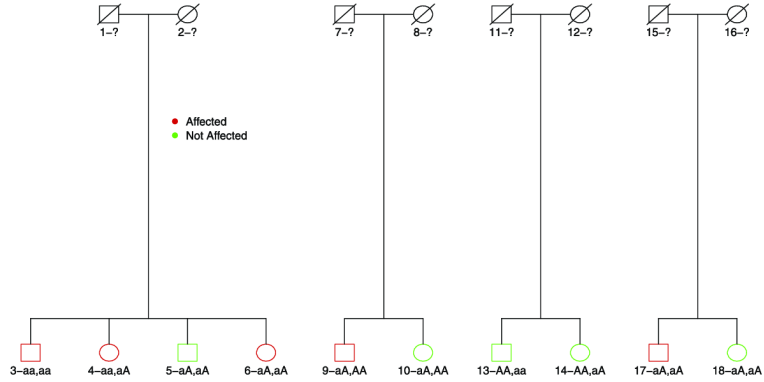}
\caption{Example of a genealogy tree with 4 families, 18 individuals
and 2 SNPs (SNP1,SNP2). Ancestors have not been observed (white); while
the offspring are labelled as healthy (green) or affected (red).}
\label{fig-exgentree}
\end{figure}

An example of such data may come from Figure \ref{fig-exgentree} in
which we have a population composed of 4 families, 18 individuals and 2
SNPs labeled as SNP1 and SNP2. From Figure \ref{fig-exgentree} we have
that ancestors have not been observed because have died (white); while
offsprings are labelled as healthy (green) or affected (red) along with
their SNP configurations.

Data in Figure \ref{fig-exgentree} may be formally represented as
follows: for individual $i$ let $Y_i \in\{0,1\}$ represent the
indicator of the phenotype, i.e. $Y_i=1$ if affected and $Y_i=0$
otherwise; assume $X_i \in\{ \{aA,Aa\},aa,AA\}$ represents the
genotype, e.g. the SNP configuration with three levels. For a genealogy
composed of $N$ individuals we have to model the corresponding pairs
$(Y_1,X_1),\ldots,(Y_N,X_N)$, where only $n<N$ have been observed. At
the phenotype level we assume the usual $\logit$ model
$(Y_i|X_i,\vtheta) \sim\mbox{Bernoulli}(p_i),$
with $p_i$ being the probability that individual $i$ is affected.

If one considers the data in Figure \ref{fig-exgentree} as $n=10$
independent observations and estimates, for instance, a logistic
regression model of $Y$ against $X$, or just considers the Fisher exact
test among $Y$ and $X$, one would end up finding no association. In
particular, the Fisher exact test between $Y$ and the first SNP has a
$p$-value of 0.21 while with SNP2 it is exactly 1. While the latter
$p$-value is reasonable as there is apparently no association between
SNP2 and $Y$, the former is not, because all $aa$ individuals are
affected and all $AA$ individuals are healthy and therefore there
should be certain evidence of association between the first SNP and
affection status. The shortcoming of this analysis is that it treats
all individuals as independent and identically distributed, while it is
clear that they are not.

On the other hand, the situation of highly dependent observations
complicates the statistical model. In fact, the probabilistic model for
the sample must take into account the genealogy that underlies the
genetic variant transmission and also the model which relates the
phenotype to the genotype. As the genotype is observed for the very
last generations only, the configurations of the SNPs for the previous
generations then constitute an enormous number of random latent
variables. This makes it almost impossible to write the likelihood for
the parameter relating the SNPs configuration with the phenotype, that
is, the coefficients of a logistic regression between $Y$ and $X$.

In this application we illustrate the \abcql\; method for the data in
Figure \ref{fig-exgentree} and we also consider a sample from the
village of Talana (Ogliastra, Sardinia - Italy) that is affected by a
reduced Mean Cell Volume (i.e. a MCV$<$72) disease for which it is
known that there exists a genetic variant inside the Beta-Globin gene
which determines it. The data, provided by the Centro Nazionale
Ricerche of Italy, consists of $N=1997$ individuals, all in one tree,
originating from two common ancestors. Only $n=49$ individuals of the
later generations are observed, among whom only 5 are
affected.
Moreover, the proportion of affected, similar to the prevalence of the
disease, is around 13\%. There are 91 SNPs with three levels and we
know that only one is inside the Beta-Globin gene.

In this analysis, we treat each SNP separately and set up a stochastic
model for a single SNP. The overall analysis for all SNPs is made by
the sequence of separated analyses over all SNPs. In the linear
predictor of the logistic regression, we consider as covariates the
genotype for individual $i$:
\[
\logit(p_i)=\theta_1 \1_{X_i=\{aA,Aa\}} +\theta_2\1_{\{X_i=aa\}
}+\theta_3 \1_{\{X_i=AA\}}.
\]

The vector of coefficients $\vtheta$ are usually interpreted as the
log of the odds ratio for the probabilities of being affected given a
SNP configuration. In order to account for the fact that the sample
$(y_1,x_1),\ldots,(y_n,x_n)$ is not i.i.d. we include the transmission
model for the genetic variants. Specifically, let $X_{i_1} $ and
$X_{i_2}$ be the SNP configuration for the ancestors of individual $i$.
Then the probabilistic model for the transmission of a genotype variant
is assumed to be regulated by the usual Mendelian inheritance model of
transmission, where the ancestor individuals are assumed to be known,
according to the genealogical tree. This law is also known as the law
of independent assortment, segregation or dominance, see for instance
\cite{levitan1988}. Therefore, if individual $i$ is a descendent in
the tree
\[
(X_i|X_{i_1},X_{i_2}) \sim\mbox{Mendel's law},
\]
while if $i$ is a founder or her/his ancestors are not in the tree,
then we assume the following prior distribution for configuration of
ancestors:
\[
X_i \sim\mbox{Trinomial}(1/3,1/3,1/3).
\]

The summary statistics, calculated only for the $n$ observations, are
the observed log-odds of the proportion of affected among all
individuals that have a certain configuration. Specifically, let $\#\{
\omega\}$ count the number of occurrences of type $\omega$,
\[
s_k=\log\left(\frac{1+\#\{ Y=1|X=k \} }{2+\#\{X=k \} }\right),
\]
where $k=1,2,3$ corresponds to $X=aa$, $X=\{aA,Aa\}$ and $X=AA$,
respectively. Note that we add one individual in the numerator and two
in the denominator in order to guarantee that $s_1,s_2$ and $s_3$ are
always defined. This, of course, constitutes a limitation for very
small samples from which, however, it would be difficult to estimate
very strong signals for a SNP that is a risk or protection factor.

In order to run \abcql, we performed a pilot-run simulation with
$M=30^3$ points on a regular grid of log odds ratios between -10 and
10. This pilot-run study depends only on the genealogy tree and not on
the observed genotypes or phenotypes. In the case of GWAS analysis,
this provides a saving in computational efforts, as in general there
are many genes to be analysed.

\begin{figure}[t]
\includegraphics{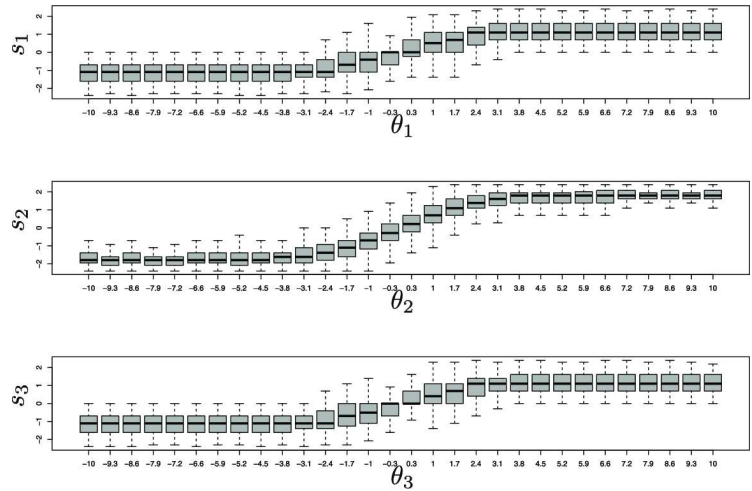}
\caption{Pilot-run for the genealogy tree in Figure \ref
{fig-exgentree} with $M=30^3$ points. Each boxplot is the conditional
distribution of $s_k|\theta_k$, $k=1,2,3$.}
\label{fig-exgentreepilotrun}
\end{figure}

The results for the pilot-run study are summarized in Figure \ref
{fig-exgentreepilotrun} where we can see that the chosen statistics are
quite informative around the null hypothesis of no association. For
very large signals, e.g. $|\theta_k|>3$, the summary statistics are
very weakly informative. This is not the fault of the summary
statistics, but it is due to the small observed sample, as is typical
in genetic isolates.

\begin{figure}
\includegraphics{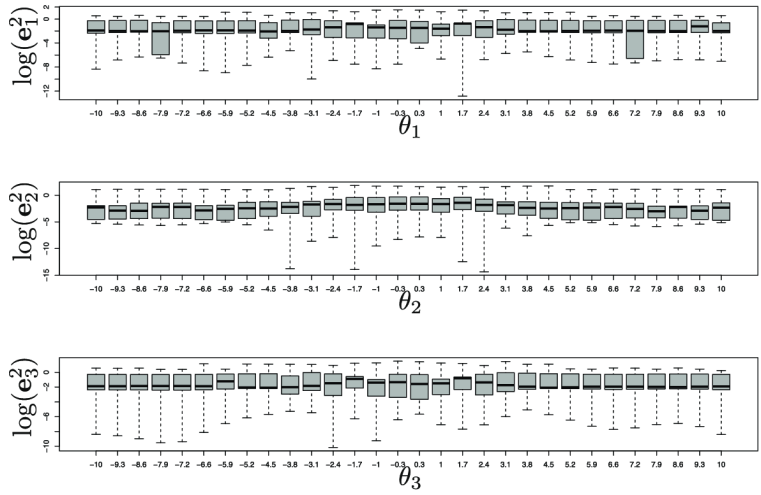}
\caption{Logarithms of squared residuals in the pilot-run for the
genealogy tree in Figure \ref{fig-exgentree}. Each boxplot is the
conditional distribution of $\log(\ve^2_k)|\theta_k$, $k=1,2,3$.
Such values indicate that $\Sr(\vtheta)$ is constant with respect to
$\vtheta$.}
\label{fig-exgentreepilotrun-var}
\end{figure}

Figure \ref{fig-exgentreepilotrun-var} illustrates the conditional
distributions of the logarithm of squared residuals with respect to
$\theta_1$, $\theta_2$ and $\theta_3$, which could be used to
estimate the variance function $\Sr(\vtheta)$. However, in this case
we find it reasonable to assume that $\Sr(\vtheta)$ is constant and
we estimate it with $\hSr$ as explained in Algorithm \ref{alg-est-fmulti}.

\begin{figure}
\includegraphics{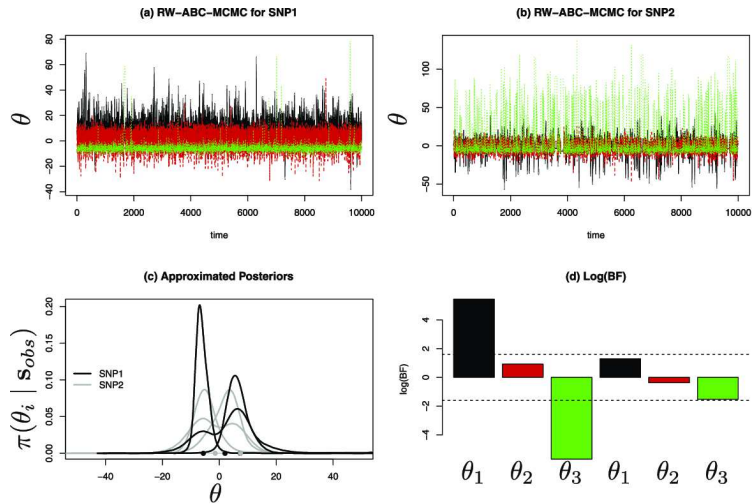}
\caption{Output of the analysis for data in Figure \ref
{fig-exgentree}. Chain output for RW-ABC-MCMC for the two SNPs (a,b)
and corresponding density estimation (c) with posterior means
represented by dots. Logarithm of Bayes Factor for $(\theta>0 \mid
\mathbf{s}_{obs})$ against $(\theta<0 \mid\mathbf{s}_{obs})$ along
with the reference lines at $\pm$1.6 (d).}
\label{fig-ex-res}
\end{figure}

We complete the algorithm by using the Euclidean distance $\rho
(\mathbf{s},\mathbf{s}_{obs})$, between $\mathbf{s}$ and its
observed value $\mathbf{s}_{obs}$ further weighted by a term that
takes into account the simulated configurations for the SNP. This term
makes the distance tend to 0 when there are many matches between
simulated and observed configurations,
\[
1-\frac{\#\{x_{obs}=x_{sim}\}}{n}.
\]
The tolerance parameter $\epsilon$ has been fixed in order to obtain
an acceptance probability in the RW-ABC-MCMC algorithm around 30\%. The
output of the chain for data in Figure \ref{fig-exgentree} is
represented in Figure \ref{fig-ex-res}.

From Figure \ref{fig-ex-res} we can see that SNP1, which has the
largest signal, exhibits values of $\theta$ with the largest posterior
mean and the largest uncertainty. Moreover, the approximated marginal
posteriors for $\theta_1$ and $\theta_3$ for SNP1 are very skewed.
For SNP2 where there is no signal, posterior distributions are centered
around 0. These results are also reflected by the logarithm of the
Bayes Factors (BFs) for $(\theta>0 \mid\mathbf{s}_{obs})$ against
$(\theta<0 \mid\mathbf{s}_{obs})$, $\Pr(\theta>0 \mid\mathbf
{s}_{obs})/\Pr(\theta<0 \mid\mathbf{s}_{obs})$ which is defined as
long as the posterior of $\theta\mid\mathbf{s}_{obs}$ exists. In
fact, there is substantial evidence for the configurations of SNP1 to
be risk or protective factors, but not for SNP2.

\begin{figure}[t!]
\includegraphics{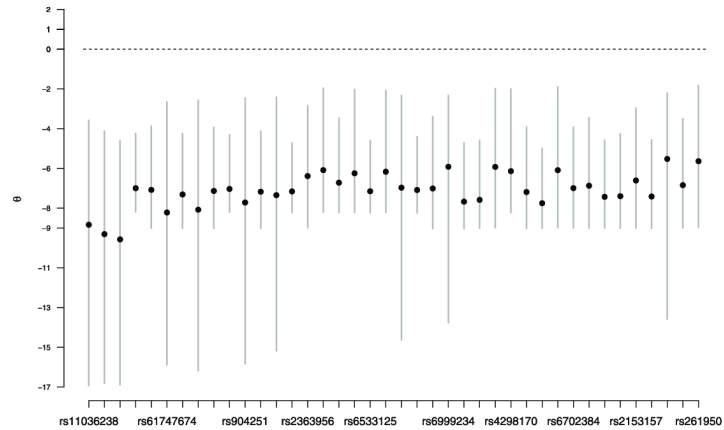}
\caption{For data from the village of Talana we report, for 40
different SNPs, 95\% credible intervals for those $\theta$s with the
largest posterior mean (in absolute) value (dot). The posterior is
approximated with the RW-ABC-MCMC.}
\label{fig-ressnp}
\end{figure}

We repeated the above analysis for the Talana data and found that the
SNP inside the Beta-Globin ({\tt rs11036238}) is among the first three
SNPs, out of 91, with the highest posterior mean (in absolute value) as
shown in Figure \ref{fig-ressnp}. Those SNPs also have the largest
Bayes Factors. However, the greater uncertainty for the first three
SNPs in Figure \ref{fig-ressnp} is due to the fact that with only
$n=49$ observed individuals we cannot be very precise in estimating
very large signals as discussed above.

\section{Conclusions}\label{sec:conclusions}
Recently, the idea of using simulation from the model to approximate
the distribution of summary statistics seems to be proliferating in the
ABC literature (see, e.g. \cite
{Prangle:2013B,Wood2010b,citeulike:12356884}). In this paper we also
used this type of approach. In particular, we simulate from the model
varying the parameters in a grid to approximate the distribution of
summary statistics as a function of the parameters in order to build a
suitable proposal for ABC-MCMC. Such a proposal distribution can be
implemented in a RW fashion or can be used as an independence kernel,
although we focused mainly on RW type MCMC algorithms.

In scalar parameter problems with conditional constant variance of
summary statistics with respect to the parameters, we showed by using
the definition of quasi-likelihood \citep{mccullagh} that this
proposal is a normal kernel in the auxiliary space, $f(\theta)$. In
multiparameter problems or when the variance of the regression function
cannot be assumed to be constant, the theory of quasi-likelihoods only
suggests a form for the proposal. In fact, analogously to the scalar
parameter case with constant variance, we propose using a multivariate
normal kernel in the auxiliary space.

A key point for the success of our method is that the summary
statistics must vary when changing the parameter values. Moreover,
there must be a one-to-one relation between them and again, the choice
of $\mathbf{s}$ is critical, as an ancillary statistic is useless for
gathering information about $\vtheta$. Such a non ancillarity
assumption is usually required over the whole parameter space $\Theta$
and it may happen that there could be parts of the parameter space
where $\mathbf{s}$ is locally ancillary. This happens, for instance,
in the coalescent model for low mutation rates and also in the
application to GWAS. Such ancillarity, however, is properly accounted
for in the discussed proposal density for ABC-MCMC. Another problem
with our approach may lie in the asymptotic argument for $p>1$ which
may not hold in some applications when $\qlik$ is irregular.

The proposed \abcql\, seems to perform quite well in the above
examples compared to other available methods; it is straightforward to
apply and its implementation does not require more than just basic
notions of regression analysis. We discussed two possible ABC-MCMC
algorithms, the RW-MH with or without constant regression variance.
Although we focus mainly on the non constant variance assumption, we
found that the original independent MH with constant regression
variance leads to satisfactory results for the discussed examples. This
is because we are only implementing a proposal distribution for the
ABC-MCMC and also because the estimation theory of QL functions allows
for a good approximation of $\nlik$ which is reflected in the proposal.

Furthermore, when simulation from $\nmodel$ is costly another
alternative to ABC could be using the $\qlik$ as a surrogate of $\nlik
$ as in \cite{a-metron-2014}. For other surrogate pseudo-likelihoods
used in the ABC context, see also \cite
{mengersen2012approximate,pauli2011composite} among others.



%

%
\begin{acknowledgement}
Maria Eugenia Castellanos was partially funded by Ministerio de Ciencia
e Innovaci\'on grant MTM2013-42323. Stefano Cabras has been partially
funded by Ministerio de Ciencia e Innovaci\'on grant MTM2013-42323,
ECO2012-38442, RYC-2012-11455 and together with Erlis Ruli were
partially funded by Ministero dell'Istruzione, dell'Univesit\`{a} e
della Ricerca of Italy. All authors have been partially financially
supported by Regione Autonoma della Sardegna under grant CRP-59903.

The authors also thank Mario Pirastu and Maria Pina Concas for
providing genealogy tree data.
\end{acknowledgement}

\end{document}